\documentclass{ifacconf}

\usepackage{graphicx}      
\usepackage{natbib}        
\usepackage{amsmath}
\usepackage{amssymb}
\usepackage{mathrsfs}
\usepackage{microtype}
\usepackage{subcaption}

\usepackage{pgfplots}
 \pgfplotsset{compat=newest}
  \usetikzlibrary{plotmarks}
  \usepackage{grffile}
\usepackage{graphics,xcolor}
\definecolor{P285U}{cmyk}{0.89,0.43,0.0,0.0}
\definecolor{P285U_font}{cmyk}{0.89,0.43,0.0,0.4}
\definecolor{lgray}{cmyk}{0,0,0,0.2}
\definecolor{myblue}{cmyk}{100,75,0,0}
\definecolor{jkuBlue}{RGB}{4,110,152}
\definecolor{jkuBlue}{RGB}{0,120,170}
\definecolor{jkuCyan}{RGB}{100,180,190}
\definecolor{jkuYellow}{RGB}{230,195,35}
\definecolor{jkuGrey}{RGB}{125,130,140}
\definecolor{jkuDarkGrey}{RGB}{51,51,51}
\definecolor{jkuLightGreen}{RGB}{195,215,75}
\definecolor{jkuGreen}{RGB}{115,180,85}
\definecolor{jkuPurple}{RGB}{145,75,130}
\definecolor{jkuRed}{RGB}{205,90,80}

\usepackage{mdframed}
\begin{document}
\begin{frontmatter}

\title{Energy-based Control and Observer Design for higher-order infinite-dimensional Port-Hamiltonian Systems\thanksref{footnoteinfo}} 

\thanks[footnoteinfo]{This work has been supported by the Austrian Science Fund (FWF) under
grant number P 29964-N32.}

\author[First]{Tobias Malzer}
\author[First]{Lukas Ecker}
\author[First]{Markus Sch\"{o}berl}

\address[First]{Institute of Automatic
Control and Control Systems Technology, Johannes Kepler University Linz,
Altenbergerstrasse 66, 4040 Linz, Austria (e-mail: \{tobias.malzer\_1, lukas.ecker, markus.schoeberl\}@jku.at).}

\begin{abstract}                
In this paper, we present a control-design method based on the energy-Casimir method for infinite-dimensional, boundary-actuated port-Hamiltonian systems with two-dimensional spatial domain and second-order Hamiltonian. The resulting control law depends on distributed system states that cannot be measured, and therefore, we additionally design an infinite-dimensional observer by exploiting the port-Hamiltonian system representation. A Kirchhoff-Love plate serves as an example in order to demonstrate the proposed approaches.  
\end{abstract}

\begin{keyword}
infinite-dimensional systems, partial differential equations, boundary actuation, port-Hamiltonian systems, structural invariants, observer design
\end{keyword}

\end{frontmatter}

\section{Introduction}

The port-Hamiltonian (pH) system representation has turned out to
be a powerful tool for the description of systems governed by ordinary
differential equations (ODEs) as well as partial differential equations
(PDEs). With respect to the infinite-dimensional case, especially
the well-known Stokes-Dirac scenario, see, e.g., \cite{Schaft2002,Gorrec2005},
is widely used, as the underlying structure --- in particular so-called
power ports --- can be exploited for the controller design like in
\cite{Macchelli2017} for instance. A key feature of this approach
is the use of energy variables replacing spatial derivatives that
occur in the Hamiltonian. This has the consequence that differential
operators, which generate the mentioned power ports, appear in the
interconnection mapping. Exemplarily, strain variables are used for
a proper pH-description of mechanical systems, which implies that
for spatially two-dimensional systems besides the PDEs also certain
compatibility conditions have to be fulfilled, see, e.g., \cite{Brugnoli2019}
for a pH-formulation of a Kirchhoff plate based on Stokes-Dirac structures.

However, from the author's point of view in particular for mechanical
systems allowing for a variational characterisation an approach based
on jet-bundle structures, see, e.g., \cite{Ennsbrunner2005,Schoeberl2015e},
is quite suitable, as the deflection of the system under consideration
appears as system state. This is especially beneficial for position
control, see, e.g., \cite{Malzer2020a} for the controller design
for infinite-dimensional systems with two-dimensional spatial domain
and in-domain actuation based on the well-known energy-Casimir method.
This pH-system formulation heavily exploits so-called jet variables
(or derivative coordinates), which are of particular importance with
respect to the generation of power ports, and therefore for the design
of boundary-control schemes, see \cite{Rams2017a} for the controller
design based on the energy-Casimir method for spatially one-dimensional
systems. In this paper, one of the intentions is to adapt the energy-Casimir
method for boundary-actuated pH-systems with $2$-dimensional spatial
domain. However, we find that for this scenario the control law depends
on system states that are distributed over a part of the boundary,
which therefore cannot be measured. In light of this aspect, a further
objective is to develop an infinite-dimensional observer by exploiting
the pH-formulation.

Note that in \cite{Toledo2019} passive observers for distributed-parameter
pH-systems are developed based on Stokes-Dirac structures. Moreover,
in \cite{Malzer2021} an observer is derived within the jet-bundle
framework for an in-domain actuated vibrating string, where the convergence
of the observer error is verified by means of functional analytic
methods. At this point, let us mention that in this contribution we
focus on energy considerations and neglect detailed stability investigations.
However, the mentioned approaches are restricted to systems with $1$-dimensional
spatial domain, whereas we present an observer-design method being
able to cope with spatially $2$-dimensional systems implying a rise
of complexity. Here, the intention is to exploit boundary-power ports
in order to impose a desired behaviour on the observer error. 

Therefore, the main contributions of this paper are as follows: i)
in Section 4, we extend the energy-Casimir method to boundary-controlled
pH-systems with $2$-dimensional spatial domain; ii) as for that scenario
the control law depends on distributed system states, an observer-design
method based on the pH-formulation is presented in Section 5. To demonstrate
the capability of the presented approaches, we study a Kirchhoff-Love
plate as running example.

\section{Notation and Preliminaries}

Throughout this contribution, differential-geometric methods are exploited
as underlying framework, where a notation similar to those of \cite{Saunders1989}
is used. Moreover, tensor notation and Einstein's convention on sums
are applied to keep formulas short and readable, where the range of
the indices is omitted when it is clear from the context. The standard
symbols $\wedge$, $\mathrm{d}$ and $\rfloor$ denote the exterior
wedge product, the exterior derivative and the Hook operator enabling
the natural contraction between tensor fields, respectively. By $C^{\infty}(\mathcal{M})$
we denote the set of all smooth functions on a manifold $\mathcal{M}$.

To properly describe distributed-parameter systems with $2$-dimensional
spatial domain in a differential-geometric setting, first we introduce
a bundle $\pi:\mathcal{E}\rightarrow\mathcal{B}$, with $(z^{1},z^{2})$
denoting the independent coordinates of the base manifold $\mathcal{B}$
and $(z^{i},x^{\alpha})$, $i=1,2$, $\alpha=1,\ldots,n$, those of
the total manifold $\mathcal{E}$. Next, we consider the $1$st-order
jet manifold $\mathcal{J}^{1}(\mathcal{E})$ equipped with the coordinates
$(z^{1},z^{2},x^{\alpha},x_{[10]}^{\alpha},x_{[01]}^{\alpha})$. Here,
we already used ordered multi-indices $[10]$ and $[01]$ representing
the $1$st-order jet variables $x_{[10]}^{\alpha}$ and $x_{[01]}^{\alpha}$,
respectively. Thus, an ordered multi index $[J]=[j_{1},j_{2}]$, with
$j_{1}+j_{2}=\#J$ and $0\leq\#J\leq r$ denoting the corresponding
order, also allows to introduce higher-order jet manifolds $\mathcal{J}^{r}(\mathcal{E})$
possessing the coordinates $(z^{1},z^{2},x_{[J]}^{\alpha})$, where
$x_{[00]}^{\alpha}=x^{\alpha}$.

A further important differential-geometric object is a tangent bundle
$\tau_{\mathcal{E}}:\mathcal{T}(\mathcal{E})\rightarrow\mathcal{E}$,
where $\mathcal{T}(\mathcal{E})$ possesses the coordinates $(z^{i},x^{\alpha},\dot{z}^{i},\dot{x}^{\alpha})$
together with the fibre basis $\partial_{i}=\partial/\partial z^{i}$,
$\partial_{\alpha}=\partial/\partial x^{\alpha}$, which allows to
introduce a vector field $v:\mathcal{E}\rightarrow\mathcal{T}(\mathcal{E})$
reading $v=v^{i}\partial_{i}+v^{\alpha}\partial_{\alpha}$ in local
coordinates. In this paper, we are particularly interested in vertical
vector fields $v:\mathcal{E}\rightarrow\mathcal{V}(\mathcal{E})$,
with $v=v^{\alpha}\partial_{\alpha}$, which can be defined by means
of vertical tangent bundles $\nu_{\mathcal{E}}:\mathcal{V}(\mathcal{E})\rightarrow\mathcal{E}$
endowed with $(z^{i},x^{\alpha},\dot{x}^{\alpha})$. Moreover, the
total derivative $d_{[1_{i}]}=\partial_{i}+x_{[J+1_{i}]}^{\alpha}\partial_{\alpha}^{[J]}$,
where $\partial_{\alpha}^{[J]}=\partial/\partial x_{[J]}^{\alpha}$
and $[1_{i}]$ represents a multi index containing only zeros except
the $i$th entry which is one, enables to prolong a vertical vector
field to the $r$th-order jet manifold according to $j^{r}(v)=v+d_{[J]}(v^{\alpha})\partial_{\alpha}^{[J]}$
with $d_{[J]}=(d_{[10]})^{j_{1}}\circ(d_{[01]})^{j_{2}}$ and $1\leq\#J\leq r$.

Next, we consider a cotangent bundle $\tau_{\mathcal{E}}^{*}:\mathcal{T}^{*}(\mathcal{E})\rightarrow\mathcal{E}$
equipped with the coordinates $(z^{i},x^{\alpha},\dot{z}_{i},\dot{x}_{\alpha})$
and the fibre basis $\mathrm{d}z^{i}$, $\mathrm{d}x^{\alpha}$, which
allows to introduce a $1$-form $\omega:\mathcal{E}\rightarrow\mathcal{T}^{*}(\mathcal{E})$
locally given as $\omega:\omega_{i}\mathrm{d}z^{i}+\omega_{\alpha}\mathrm{d}x^{\alpha}$.
Moreover, in this paper we study Hamiltonian densities $\mathfrak{H}=\mathcal{H}\Omega$,
which depend on $2$nd-order jet variables, i.e. $\mathcal{H}\in C^{\infty}(\mathcal{J}^{2}(\mathcal{E}))$,
and can be constructed by means of special pullback bundles omitted
here for ease of presentation. Here, $\Omega=\mathrm{d}z^{1}\wedge\mathrm{d}z^{2}$
denotes a volume form, whereas $\Omega_{i}=\partial_{i}\rfloor\Omega$
corresponds to a boundary-volume form. Moreover, the formal change
of a Hamiltonian functional $\mathscr{H}=\int_{\mathcal{B}}\mathcal{H}\Omega$
along the solutions of an evolutionary vector field $v=v^{\alpha}\partial_{\alpha}$,
corresponding to a set of PDEs $\dot{x}^{\alpha}=v^{\alpha}$ with
$v^{\alpha}\in C^{\infty}(\mathcal{J}^{4}(\mathcal{E}))$ and the
time $t$ as evolution parameter of the solution, is of great significance,
where we use the Lie-derivative reading $\mathrm{L}_{v}(\omega)$
for a differential form $\omega$. Due to the fact that $\mathcal{H}\in C^{\infty}(\mathcal{J}^{2}(\mathcal{E}))$,
we are interested in $\dot{\mathscr{H}}=\int_{\mathcal{B}}\mathrm{L}_{j^{2}v}(\mathcal{H}\Omega)$.
In fact, for the system configuration under investigation, i.e. $2$nd-order
Hamiltonian density and $2$-dimensional spatial domain, the determination
of $\dot{\mathscr{H}}$ is a non-trivial task. However, in \cite{Schoeberl2018}
an approach based on so-called Cartan forms is presented allowing
to introduce the boundary operators\begin{subequations}\label{eq:boundary_operators}
\begin{align}
\delta^{\partial,1}\mathfrak{H} & =(\partial_{\alpha}^{[01]}\mathcal{H}-d_{[10]}(\partial_{\alpha}^{[11]}\mathcal{H})-d_{[01]}(\partial_{\alpha}^{[02]}\mathcal{H}))\mathrm{d}x^{\alpha}\wedge\bar{\Omega}_{2}\\
\delta^{\partial,2}\mathfrak{H} & =\partial_{\alpha}^{[02]}\mathcal{H}\mathrm{d}x_{[01]}^{\alpha}\wedge\bar{\Omega}_{2},
\end{align}
\end{subequations}where $\bar{\Omega}_{2}$ denotes a boundary-volume
form in coordinates adapted to the boundary. Moreover, for $\mathcal{H}\in C^{\infty}(\mathcal{J}^{2}(\mathcal{E}))$
and $\mathrm{dim}(\mathcal{B})=2$, the variational derivative reads
\begin{equation}
\delta\mathfrak{H}=\delta_{\alpha}\mathcal{H}\mathrm{d}x^{\alpha}\wedge\Omega\label{eq:domain_operator}
\end{equation}
with $\delta_{\alpha}=\partial_{\alpha}-d_{[10]}\partial_{\alpha}^{[10]}-d_{[01]}\partial_{\alpha}^{[01]}+d_{[20]}\partial_{\alpha}^{[20]}+d_{[11]}\partial_{\alpha}^{[11]}+d_{[02]}\partial_{\alpha}^{[02]}$.
Thus, based on \cite{Schoeberl2018}, it is possible to introduce
the following theorem.

\begin{thm}[Decomposition Theorem]\label{thm:decomposition_theorem}\cite[Theorem 3.2]{Rams2018}
Let $v$ be an evolutionary vector field and $\mathcal{H}\in C^{\infty}(\mathcal{J}^{2}(\mathcal{E}))$
a second-order density. Then, by exploiting the domain operator (\ref{eq:domain_operator})
and the boundary operators (\ref{eq:boundary_operators}), the integral
$\int_{\mathcal{B}}\mathrm{L}_{j^{2}v}(\mathcal{H}\Omega)$ can be
decomposed into
\[
\dot{\mathscr{H}}=\int_{\mathcal{B}}v\rfloor\delta\mathfrak{H}+\int_{\partial\mathcal{B}}v\rfloor\delta^{\partial,1}\mathfrak{H}+\int_{\partial\mathcal{B}}j^{1}(v)\rfloor\delta^{\partial,2}\mathfrak{H}.
\]
\end{thm}

Finally, let us mention that the bundle $\pi:\mathcal{E}\rightarrow\mathcal{B}$
enables the construction of further differential-geometric objects
such as the tensor bundle $\mathcal{W}_{1}^{2,r}=\mathcal{T}^{*}(\mathcal{E})\wedge\bigwedge^{2}\mathcal{T}^{*}(\mathcal{B})$,
where $\omega=\omega_{\alpha}\mathrm{d}x^{\alpha}\wedge\Omega$, with
$\omega_{\alpha}\in C^{\infty}(\mathcal{J}^{r}(\mathcal{E}))$, denotes
an element.

\section{Infinite-Dimensional PH-Systems}

Next, we discuss the port-Hamiltonian framework for systems with $2$nd-order
Hamiltonian, see \cite{Rams2017a}, where we focus on systems with
$2$-dimensional spatial domain. The presented approach is based on
an underlying jet-bundle structure as well as on a certain power-balance
relation, where we exploit Theorem \ref{thm:decomposition_theorem}.

\begin{figure}
\centering
\def \svgwidth{0.45\textwidth}
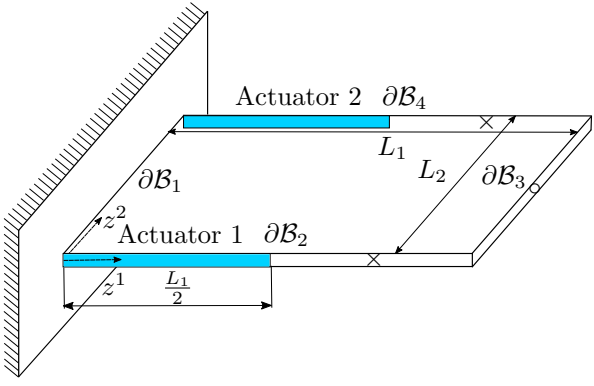\caption{\label{fig:Schematic_Kirchhoff_Love_plate}Schematic representation
of the boundary-actuated Kirchhoff-Love plate.}
\end{figure}

Thus, we consider systems with $2$nd-order Hamiltonian $\mathfrak{H}$,
i.e. $\mathcal{\ensuremath{H}}\in C^{\infty}(\mathcal{J}^{2}(\mathcal{E}))$,
and a $2$-dimensional, rectangular spatial domain $\mathcal{B}=\{(z^{1},z^{2})|z^{1}\in[0,L_{1}],z^{2}\in[0,L_{2}]\}$,
where the boundary is divided into $\partial\mathcal{B}_{1}=\{(z^{1},z^{2})|z^{1}=0,z^{2}\in[0,L_{2}]\}$,
$\partial\mathcal{B}_{2}=\{(z^{1},z^{2})|z^{1}\in[0,L_{1}],z^{2}=0\}$,
$\partial\mathcal{B}_{3}=\{(z^{1},z^{2})|z^{1}=L_{1},z^{2}\in[0,L_{2}]\}$
and $\partial\mathcal{B}_{4}=\{(z^{1},z^{2})|z^{1}\in[0,L_{1}],z^{2}=L_{2}\}$,
see Fig. \ref{fig:Schematic_Kirchhoff_Love_plate}. Then, a pH-formulation
can be given as
\begin{equation}
\dot{x}=(\mathcal{J}-\mathcal{R})(\delta\mathfrak{H})\label{eq:boundary_controlled_ipH_non_diff_op}
\end{equation}
together with appropriate boundary conditions. The objects $\mathcal{J}$
and $\mathcal{R}$ describe the internal power flow and dissipation
effects, respectively, and can be interpreted as mappings of the form
$\mathcal{J},\mathcal{R}:\mathcal{W}_{1}^{2,4}(\mathcal{E})\rightarrow\mathcal{V}(\mathcal{E})$.
Next, the formal change of the Hamiltonian functional $\mathscr{H}$
along solutions of (\ref{eq:boundary_controlled_ipH_non_diff_op}),
which can be written as
\begin{equation}
\dot{\mathscr{H}}=-\int_{\mathcal{B}}\mathcal{R}(\delta\mathfrak{H})\rfloor\delta\mathfrak{H}+\int_{\partial\mathcal{B}}(\dot{x}\rfloor\delta^{\partial,1}\mathfrak{H}+\dot{x}_{[01]}\rfloor\delta^{\partial,2}\mathfrak{H})\label{eq:formal_change_non_diff_op}
\end{equation}
by substituting $v=\dot{x}$ with (\ref{eq:boundary_controlled_ipH_non_diff_op})
in Theorem \ref{thm:decomposition_theorem}, is of particular interest
and states a power-balance relation if $\mathscr{H}$ corresponds
to the total energy of the system. Moreover, a local coordinate representation
of (\ref{eq:boundary_controlled_ipH_non_diff_op}) can be given by
\begin{equation}
\dot{x}^{\alpha}=(\mathcal{J}^{\alpha\beta}-\mathcal{R}^{\alpha\beta})\delta_{\beta}\mathcal{H},\qquad\alpha,\beta=1,\ldots,n.\label{eq:boundary_controlled_ipH_non_diff_op_local}
\end{equation}
While the interconnection map $\mathcal{J}$ is skew-symmetric, i.e.
the coefficients fulfil $\mathcal{J}^{\alpha\beta}=-\mathcal{J}^{\beta\alpha}\in C^{\infty}(\mathcal{J}^{4}(\mathcal{E}))$,
the dissipation map $\mathcal{R}$ is symmetric and positive semi-definite,
implying that the coefficients meet $\mathcal{R}^{\alpha\beta}=\mathcal{R}^{\beta\alpha}\in C^{\infty}(\mathcal{J}^{4}(\mathcal{E}))$
and $[\mathcal{R}^{\alpha\beta}]\geq0$ for the coefficient matrix.
Thus, the power-balance relation
\[
\dot{\mathscr{H}}\!=\!-\!\int_{\mathcal{B}}\delta_{\alpha}(\mathcal{H})\mathcal{R}^{\alpha\beta}\delta_{\beta}(\mathcal{H})\Omega+\int_{\partial\mathcal{B}}(\dot{x}^{\alpha}\delta_{\alpha}^{\partial,1}\mathcal{H}+\dot{x}_{[01]}^{\alpha}\delta_{\alpha}^{\partial,2}\mathcal{H})\bar{\Omega}_{2}
\]
allows to introduce boundary ports in a straightforward manner by
exploiting the boundary operators (\ref{eq:boundary_operators}).
Now, we assume that the boundary $\partial\mathcal{B}$ is divided
into an actuated boundary $\partial\mathcal{B}_{a}=\partial\mathcal{B}_{2}\cup\partial\mathcal{B}_{4}$
and an unactuated boundary $\partial\mathcal{B}_{u}=\partial\mathcal{B}_{1}\cup\partial\mathcal{B}_{3}$.
This has the consequence that for the actuated boundary $\partial\mathcal{B}_{a}$
we are able to set
\[
(\dot{x}^{\alpha}\delta_{\alpha}^{\partial,1}\mathcal{H})|_{\partial\mathcal{B}_{a}}=u_{\partial,1}^{\varsigma}y_{\varsigma}^{\partial,1},\quad(\dot{x}_{[01]}^{\alpha}\delta_{\alpha}^{\partial,2}\mathcal{H})|_{\partial\mathcal{B}_{a}}=u_{\partial,2}^{\mu}y_{\mu}^{\partial,2},
\]
while no power flow takes place over the unactuated boundary $\partial\mathcal{B}_{u}$,
i.e. $(\dot{x}^{\alpha}\delta_{\alpha}^{\partial,1}\mathcal{H})|_{\partial\mathcal{B}_{u}}=0$
and $(\dot{x}_{[01]}^{\alpha}\delta_{\alpha}^{\partial,2}\mathcal{H})|_{\partial\mathcal{B}_{u}}=0$.
It should be noted that the roles of inputs and outputs cannot be
uniquely defined; however, for our purposes, where we intend to exploit
mechanical quantities such as forces and bending moments as inputs
and velocities and angular velocities, respectively, as collocated
output quantities, we confine ourselves to the parameterisation\begin{subequations}\label{eq:ipH_sys_boundary_parameterisation}
\begin{equation}
\begin{array}{rclcrcl}
\delta_{\alpha}^{\partial,1}\mathcal{H}|_{\partial\mathcal{B}_{a}} & = & B_{\alpha\varsigma}^{\partial,1}u_{\partial,1}^{\varsigma}, & \quad & \delta_{\alpha}^{\partial,2}\mathcal{H}|_{\partial\mathcal{B}_{a}} & = & B_{\alpha\mu}^{\partial,2}u_{\partial,2}^{\mu},\end{array}\label{eq:ipH_sys_input_parameterisation}
\end{equation}
for the boundary inputs as well as
\begin{equation}
\begin{array}{ccccccc}
B_{\alpha\varsigma}^{\partial,1}\dot{x}^{\alpha}|_{\partial\mathcal{B}_{a}} & = & y_{\varsigma}^{\partial,1} & \quad & B_{\alpha\mu}^{\partial,2}\dot{x}_{[01]}^{\alpha}|_{\partial\mathcal{B}_{a}} & = & y_{\mu}^{\partial,2},\end{array}\label{eq:ipH_sys_output_parameterisation}
\end{equation}
\end{subequations}for the outputs, where $\varsigma=1,\ldots,l_{\partial,1}$,
$\mu=1,\ldots,l_{\partial,2}$.
\begin{figure}
\centering
\input{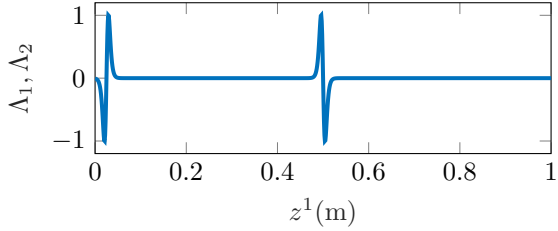}\caption{\label{fig:Force_distribution_actuators}Force distribution of the
actuators over $z^{1}$.}
\end{figure}

\begin{exmp}[Boundary-actuated Kirchhoff-Love plate]\label{ex:Kirchhoff_Love_plate}We
\newline consider a rectangular, boundary-actuated Kirchhoff-Love
plate depicted in Fig. \ref{fig:Schematic_Kirchhoff_Love_plate} and
governed by
\begin{equation}
\rho A\ddot{w}=-D_{E}(w_{[40]}+2w_{[22]}+w_{[04]})\label{eq:Kirchhoff_Plate_PDE}
\end{equation}
see \cite[p. 448, Eq. (7.333)]{Meirovitch1997}. The boundary conditions
are discussed below, as they are of particular interest here. For
the sake of simplicity we assume the material parameters $\rho,A,D_{E}>0$
in (\ref{eq:Kirchhoff_Plate_PDE}) to be constant. If we introduce
the generalised momentum $p=\rho A\dot{w}$ as well as the Hamiltonian
density $\mathcal{H}=\mathcal{T}+\mathcal{V}$, where $\mathcal{T}=\frac{1}{2\rho A}p^{2}$,
\[
\mathcal{V}\!=\!\tfrac{1}{2}D_{E}((w_{\left[20\right]})^{2}+(w_{\left[02\right]})^{2}+2\nu w_{\left[20\right]}w_{\left[02\right]}+2(1-\nu)(w_{\left[11\right]})^{2})
\]
with Poisson's ratio $\nu$, we find an appropriate pH-system representation
according to
\[
\left[\begin{array}{c}
\dot{w}\\
\dot{p}
\end{array}\right]\!=\!\left[\begin{array}{cc}
0 & 1\\
-1 & 0
\end{array}\right]\left[\begin{array}{c}
\delta_{w}\mathcal{H}\\
\delta_{p}\mathcal{H}
\end{array}\right]\!=\!\left[\begin{array}{c}
\tfrac{p}{\rho A}\\
-D_{E}(w_{[40]}+2w_{[22]}+w_{[04]})
\end{array}\right].
\]
Thus, the power-balance relation (\ref{eq:formal_change_non_diff_op})
reads
\begin{multline}
\dot{\mathscr{H}}=\int_{0}^{L_{1}}(Q_{1}\dot{w}+M_{1}\dot{w}_{[01]})|_{z^{2}\in\{0,L_{2}\}}\mathrm{d}z^{1}\\
+\int_{0}^{L_{2}}(Q_{2}\dot{w}+M_{2}\dot{w}_{[10]})|_{z^{1}\in\{0,L_{1}\}}\mathrm{d}z^{2},\label{eq:Kelvin_Voigt_plate_formal_change_1}
\end{multline}
where the boundary relations $Q_{1}=D_{E}(w_{[03]}+(2-\nu)w_{[21]})$,
$M_{1}=-D_{E}(w_{[02]}+\nu w_{[20]})$, $Q_{2}=-D_{E}(w_{[30]}+(2-\nu)w_{[12]})$
and $M_{2}=D_{E}(w_{[20]}+\nu w_{[02]})$ can be determined by evaluating
the boundary operators (\ref{eq:boundary_operators})\footnote{Note that for the boundaries $\partial\mathcal{B}_{2}$ and $\partial\mathcal{B}_{4}$
we have the adapted coordinates $z_{\partial}^{1}=z^{1}$ and $z_{\partial}^{2}=z^{2}=\text{const.}$,
i.e. the boundary operators indeed read (\ref{eq:boundary_operators});
however for $\partial\mathcal{B}_{1}$ and $\partial\mathcal{B}_{3}$
the coordinates adapted to the boundary are $z_{\partial}^{1}=z^{2}$
and $z_{\partial}^{2}=z^{1}=\text{const.}$ implying that the coordinates
in (\ref{eq:boundary_operators}) need to be swapped.}. 

With regard to control-engineering purposes, the boundary ports in
(\ref{eq:Kelvin_Voigt_plate_formal_change_1}) can be used to extract
or deliver power. In fact, we assume that the plate is clamped at
the boundary $\partial\mathcal{B}_{1}$ and free at $\partial\mathcal{B}_{3}$,
i.e. we have the boundary conditions
\begin{equation}
\begin{array}{rccclcccccc}
\left.\begin{array}{rcc}
w & = & 0\\
w_{[10]} & = & 0
\end{array}\right\}  & \; & \text{for} & \; & \partial\mathcal{B}_{1}, & \; & \left.\begin{array}{rcc}
Q_{2} & = & 0\\
M_{2} & = & 0
\end{array}\right\}  & \; & \text{for} & \; & \partial\mathcal{B}_{3},\end{array}\label{eq:Kirchhoff_plate_boundary conditions}
\end{equation}
implying that the corresponding power ports vanish identically. However,
the shear forces $Q_{1}|_{\partial\mathcal{B}_{2}}$ and $Q_{1}|_{\partial\mathcal{B}_{4}}$
at the actuated boundary shall be generated by piezo-like actuators,
which are supposed to be perfectly attached at the boundaries $\partial\mathcal{B}_{2}$
and $\partial\mathcal{B}_{4}$. The forces supplied by these actuators,
which are spatially distributed over (a part of) $\partial\mathcal{B}_{2}$
and $\partial\mathcal{B}_{4}$, see Fig. \ref{fig:Schematic_Kirchhoff_Love_plate}
and \ref{fig:Force_distribution_actuators}, can be described by $Q_{1}|_{\partial\mathcal{B}_{2}}=\varLambda_{1}u_{in}^{1}$
and $Q_{1}|_{\partial\mathcal{B}_{4}}=\varLambda_{2}u_{in}^{2}$,
where the characteristic functions read
\begin{align*}
\varLambda_{1} & =-\Psi\tfrac{\partial^{2}}{\partial(z^{1})^{2}}(\tanh(\sigma z^{1})-\tanh(\sigma(z^{1}-\tfrac{L_{1}}{2}))|_{z^{2}=0}\\
\varLambda_{2} & =-\Psi\tfrac{\partial^{2}}{\partial(z^{1})^{2}}(\tanh(\sigma z^{1})-\tanh(\sigma(z^{1}-\tfrac{L_{1}}{2}))|_{z^{2}=L_{2}}
\end{align*}
with the material parameters hidden in $\Psi$ and $\sigma\in\mathbb{R}^{+}$.
Thus, the voltages $u_{in}^{1}$, $u_{in}^{2}$ applied to the piezo-like
actuators serve as manipulated variables. In light of this aspects,
for (\ref{eq:ipH_sys_boundary_parameterisation}) we set
\begin{equation}
\begin{array}{ccccccc}
Q_{1}|_{\partial\mathcal{B}_{2}} & = & B_{11}^{\partial,1}u_{\partial,1}^{1}, & \qquad & B_{11}^{\partial,1}\dot{w}|_{\partial\mathcal{B}_{2}} & = & y_{1}^{\partial,1},\\
Q_{1}|_{\partial\mathcal{B}_{4}} & = & B_{12}^{\partial,1}u_{\partial,1}^{2}, & \qquad & B_{12}^{\partial,1}\dot{w}|_{\partial\mathcal{B}_{4}} & = & y_{2}^{\partial,1},
\end{array}\label{eq:param_input_output_Kirchhoff_plate_Kelvin_Voigt}
\end{equation}
with $B_{11}^{\partial,1}=\varLambda_{1}$, $B_{12}^{\partial,1}=\varLambda_{2}$,
 and $u_{\partial,1}^{1}=u_{in}^{1}$, $u_{\partial,1}^{2}=u_{in}^{2}$,
to assign the roles of inputs and outputs. Note that since we have
the shear forces $Q_{1}|_{\partial\mathcal{B}_{2}}$, $Q_{1}|_{\partial\mathcal{B}_{4}}$
as inputs solely, the relations $M_{1}|_{\partial\mathcal{B}_{2}}=0$,
$M_{1}|_{\partial\mathcal{B}_{4}}=0$ complete the boundary conditions
(\ref{eq:Kirchhoff_plate_boundary conditions}). Hence, by virtue
of the plate configuration and the assignment (\ref{eq:param_input_output_Kirchhoff_plate_Kelvin_Voigt}),
the power-balance relation reads
\[
\dot{\mathscr{H}}=\int_{\partial\mathcal{B}_{2}}u_{\partial,1}^{1}y_{1}^{\partial,1}\mathrm{d}z^{1}+\int_{\partial\mathcal{B}_{4}}u_{\partial,1}^{2}y_{2}^{\partial,1}\mathrm{d}z^{1}.
\]
\end{exmp}

\section{Boundary Control based on Structural Invariants}

Now, the objective is to adapt the energy-Casimir method to boundary-actuated
pH-systems with $2$nd-order Hamiltonian and $2$-dimensional spatial
domain, where the intention is to design a dynamic controller for
the boundary-actuated Kirchhoff-Love plate discussed in Ex. \ref{ex:Kirchhoff_Love_plate}
by exploiting a certain interconnection of plant and controller.

\subsection{Control by Interconnection}

Thus, we develop a dynamic controller based on structural invariants,
which allows to shape the total energy of the closed loop and to inject
additional damping in order to increase the dissipation rate. In light
of the aspect that we only consider systems with lumped inputs, the
interconnection of a plant (\ref{eq:boundary_controlled_ipH_non_diff_op})
and a finite-dimensional controller, beneficially given in the pH-formulation\begin{subequations}\label{eq:fpH_controller}
\begin{align}
\dot{x}_{c}^{\alpha_{c}} & =(J_{c}^{\alpha_{c}\beta_{c}}\!-\!R_{c}^{\alpha_{c}\beta_{c}})\partial_{\beta_{c}}H_{c}\!+\!G_{c,\varsigma}^{\partial,1,\alpha_{c}}u_{c,\partial,1}^{\varsigma}\!+\!G_{c,\mu}^{\partial,2,\alpha_{c}}u_{c,\partial,2}^{\mu}\label{eq:fpH_controller_dynamics}\\
y_{c,\varsigma}^{\partial,1} & =G_{c,\varsigma}^{\partial,1,\alpha_{c}}\partial_{\alpha_{c}}H_{c},\quad y_{c,\mu}^{\partial,2}=G_{c,\mu}^{\partial,2,\alpha_{c}}\partial_{\alpha_{c}}H_{c},
\end{align}
\end{subequations}with $\alpha_{c},\beta_{c}=1,\ldots,n_{c}$, $\varsigma=1,\ldots,l_{\partial,1}$
and $\mu=1,\ldots,l_{\partial,2}$, is motivated, where we splitted
the controller inputs into two different parts to take into account
the two different categories of boundary ports. The idea is to couple
the finite-dimensional controller to the infinite-dimensional plant
at the actuated boundary $\partial\mathcal{B}_{a}=\partial\mathcal{B}_{2}\cup\partial\mathcal{B}_{4}$
in a power-conserving manner fulfilling
\begin{equation}
\int_{\partial\mathcal{B}_{a}}(u_{\partial,1}\rfloor y^{\partial,1}+u_{\partial,2}\rfloor y^{\partial,2})+u_{c,\partial,1}\rfloor y_{c}^{\partial,1}+u_{c,\partial,2}\rfloor y_{c}^{\partial,2}=0.\label{eq:PCIS_2D}
\end{equation}
Note that this is quite different compared to \cite{Rams2017a}, where
boundary-actuated systems with $1$-dimensional spatial domain are
considered, as well as to \cite{Malzer2020a}, where the energy-Casimir
method is investigated for in-domain actuated systems with $2$-dimensional
spatial domain, as we have to integrate over the $1$-dimensional
boundary here. Hence, if we choose
\begin{equation}
\begin{array}{ccccccc}
u_{c,\partial,1} & = & \int_{\partial\mathcal{B}_{a}}K_{\partial,1}\rfloor y^{\partial,1}, & \quad & u_{\partial,1} & = & -K_{\partial,1}^{*}\rfloor y_{c}^{\partial,1}\\
u_{c,\partial,2} & = & \int_{\partial\mathcal{B}_{a}}K_{\partial,2}\rfloor y^{\partial,2}, & \quad & u_{\partial,2} & = & -K_{\partial,2}^{*}\rfloor y_{c}^{\partial,2},
\end{array}\label{eq:PCIS_feedback_2D}
\end{equation}
where $K_{\partial,1}$ and $K_{\partial,2}$ denote appropriate mappings
-- that can be interpreted as degrees of freedom but are set to the
identity matrix for the most part -- the relation (\ref{eq:PCIS_2D})
is satisfied. As a consequence, the closed-loop system can be formulated
as a pH-system described by the Hamiltonian $\mathscr{H}_{cl}=\mathscr{H}+H_{c}$.
Moreover, because of the power-conserving interconnection (\ref{eq:PCIS_2D}),
the formal change of $\mathscr{H}_{cl}$ along solutions of the closed-loop
system reads
\[
\dot{\mathscr{H}}_{cl}=-\int_{\mathcal{B}}\delta_{\alpha}(\mathcal{H})\mathcal{R}^{\alpha\beta}\delta_{\beta}(\mathcal{H})\Omega-\partial_{\alpha_{c}}(H_{c})R_{c}^{\alpha_{c}\beta_{c}}\partial_{\beta_{c}}(H_{c}),
\]
highlighting that the controller allows to inject damping.

\begin{rem}\label{rem:stability}At this point let us stress that
we assume that the closed-loop solutions exist. In principle, the
well-posedness of the closed-loop system would need to be verified
by means of functional analysis. However, in this paper the emphasis
is on a formal, geometric approach focusing on energy considerations,
where the relations $\mathscr{H}_{cl}>0$ and $\dot{\mathscr{H}}_{cl}\leq0$
serve as necessary conditions for (possible) stability investigations
in the sense of Liapunov. \end{rem}

\subsection{Controller Design}

Next, we are interested in certain functionals allowing for a relation
between plant and controller states in order to shape the Hamiltonian
of the closed loop. Thus, in accordance with \cite{Rams2017a}, we
consider
\begin{equation}
\mathscr{C}^{\lambda}=x_{c}^{\lambda}+\int_{\mathcal{B}}\mathcal{C}^{\lambda}\Omega,\qquad\mathcal{C}^{\lambda}\in C^{\infty}(\mathcal{J}^{2}(\mathcal{E})),\label{eq:Casimir_functions_fpH_controller_2D_2D_plant}
\end{equation}
with $\lambda=1,\ldots,\bar{n}_{c}\leq n_{c}$, where it should be
mentioned that here we have $\text{dim}(\mathcal{B})=2$. Thus, the
functionals (\ref{eq:Casimir_functions_fpH_controller_2D_2D_plant})
have to meet $\dot{\mathscr{C}}^{\lambda}=0$ independently of $\mathcal{H}$
and $H_{c}$ to qualify as structural invariants.\begin{thm}[Structural Invariants]\label{thm:Casimir_conditions}Consider
the closed-loop system stemming from the interconnection of the plant
(\ref{eq:boundary_controlled_ipH_non_diff_op}) and the controller
(\ref{eq:fpH_controller}) by means of (\ref{eq:PCIS_feedback_2D}).
Thus, (\ref{eq:Casimir_functions_fpH_controller_2D_2D_plant}) are
structural invariants if they meet the conditions\begin{subequations}\label{eq:Casimir_conditions_boundary_act_diff_op_scen}
\begin{align}
(J_{c}^{\lambda\beta_{c}}-R_{c}^{\lambda\beta_{c}}) & =0\\
\delta_{\alpha}\mathcal{C}^{\lambda}(\mathcal{J}^{\alpha\beta}-\mathcal{R}^{\alpha\beta}) & =0\label{eq:Casimir_cond_BA_diff_op_domain}\\
(G_{c,\varsigma}^{\partial,1,\lambda}K_{\partial,1}^{\varsigma\eta}B_{\alpha\eta}^{\partial,1}+\delta_{\alpha}^{\partial,1}\mathcal{C}^{\lambda})|_{\partial\mathcal{B}_{a}} & =0\label{eq:Casimir_cond_BA_diff_op_act_bound_1}\\
(G_{c,\mu}^{\partial,2,\lambda}K_{\partial,2}^{\mu\kappa}B_{\alpha\kappa}^{\partial,2}+\delta_{\alpha}^{\partial,2}\mathcal{C}^{\lambda})|_{\partial\mathcal{B}_{a}} & =0\label{eq:Casimir_cond_BA_diff_op_act_bound_2}\\
(\dot{x}^{\alpha}\delta_{\alpha}^{\partial,1}\mathcal{C}^{\lambda}+\dot{x}_{[01]}^{\alpha}\delta_{\alpha}^{\partial,2}\mathcal{C}^{\lambda})|_{\partial\mathcal{B}_{u}} & =0.\label{eq:Casimir_cond_BA_diff_op_unact_bound}
\end{align}
\end{subequations}\end{thm}

\begin{pf}To prove the conditions (\ref{eq:Casimir_conditions_boundary_act_diff_op_scen}),
we substitute $v$ and $\mathfrak{H}$ by $\dot{x}$ and $\mathcal{C}^{\lambda}\Omega$,
respectively, in the decomposition Theorem \ref{thm:decomposition_theorem}.
Moreover, by inserting the plant and controller dynamics described
by (\ref{eq:boundary_controlled_ipH_non_diff_op_local}) and (\ref{eq:fpH_controller_dynamics}),
respectively, as well as the coupling (\ref{eq:PCIS_feedback_2D})
together with the boundary-output assignments (\ref{eq:ipH_sys_output_parameterisation}),
we are able to deduce
\begin{multline*}
\dot{\mathscr{C}}^{\lambda}=(J_{c}^{\lambda\beta_{c}}-R_{c}^{\lambda\beta_{c}})\partial_{\beta_{c}}H_{c}+\int_{\partial\mathcal{B}_{a}}G_{c,\varsigma}^{\partial,1,\lambda}K_{\partial,1}^{\varsigma\eta}B_{\alpha\eta}^{\partial,1}\dot{x}^{\alpha}\bar{\Omega}_{2}\\
+\int_{\partial\mathcal{B}_{a}}G_{c,\mu}^{\partial,2,\lambda}K_{\partial,2}^{\mu\kappa}B_{\alpha\kappa}^{\partial,2}\dot{x}_{[01]}^{\alpha}\bar{\Omega}_{2}+\int_{\mathcal{B}}\delta_{\alpha}(\mathcal{C}^{\lambda})(\mathcal{J}^{\alpha\beta}-\mathcal{R}^{\alpha\beta})\delta_{\beta}\mathcal{H}\Omega\\
+\int_{\partial\mathcal{B}}\dot{x}^{\alpha}\delta_{\alpha}^{\partial,1}(\mathcal{C}^{\lambda})\bar{\Omega}_{2}+\int_{\partial\mathcal{B}}\dot{x}_{[01]}^{\alpha}\delta_{\alpha}^{\partial,1}(\mathcal{C}^{\lambda})\bar{\Omega}_{2}=0,
\end{multline*}
enabling to find the conditions (\ref{eq:Casimir_conditions_boundary_act_diff_op_scen}).\end{pf}

Note that -- in contrast to \cite[Eq. (17)]{Rams2017a}, where pH-systems
with $1$-dimensional spatial domain are considered -- the conditions
(\ref{eq:Casimir_conditions_boundary_act_diff_op_scen}) basically
hold for systems with $1$- or $2$-dimensional spatial domain; however,
the differences are hidden in the geometric objects that of course
strongly depend on the spatial dimension. Moreover, the conditions
(\ref{eq:Casimir_conditions_boundary_act_diff_op_scen}) clearly distinguish
from those of Prop. 2 in \cite{Malzer2020a}, where Casimir conditions
for in-domain actuated systems with $2$-dimensional spatial domain
are studied.

\begin{rem}It should be stressed that the proposed control scheme
can be exploited for nonlinear systems as well, see \cite[Sec. 4]{Malzer2018},
where a Casimir-based controller for a nonlinear Euler-Bernoulli beam
structure is presented. However, in this paper we intend to combine
the proposed controller with an infinite-dimensional observer, where
the design method is restricted to the linear scenario, and therefore,
in the following we derive a controller for a linear Kirchhoff-Love
plate.\end{rem}

\begin{exmp}[Casimir-based Controller for Ex. 2]\label{Ex:Casimir_Controller_Kirchhoff_Love_plate}Next,
\newline we design a controller for the Kirchhoff-Love plate to stabilise
the (approximated)\footnote{Due to the special force characteristic, our actuators are able to
generate equivalent bending moments.} configuration
\begin{equation}
w^{d}=\left\{ \begin{array}{lcccl}
a(z^{1})^{2} & \; & \text{for} & \; & 0\leq z^{1}<\tfrac{L_{1}}{2}\\
b(z^{1}-\tfrac{L_{1}}{2})+a\tfrac{L_{1}^{2}}{4} & \; & \text{for} & \; & \tfrac{L_{1}}{2}\leq z^{1}\leq L_{1}
\end{array}\right.\label{eq:BA_Kirchhoff_Plate_des_equilibrium}
\end{equation}
with $a,b>0$. In light of the fact that for the system under investigation
we have two output densities given in (\ref{eq:param_input_output_Kirchhoff_plate_Kelvin_Voigt}),
we intend to relate two controller states to the plant. To this end,
we make the trivial choice $K_{\partial,1}=I$, with $I$ denoting
the identity matrix, and $K_{\partial,2}=0$ (as we have no boundary-actuation
corresponding to the category $\partial,2$) regarding the design
parameters, and thus, the interconnection of plant and controller
reads\begin{subequations}\label{eq:BA_Kirchhoff_Plate_interconnection}
\begin{align}
u_{c,\partial,1}^{1} & =\int_{\partial\mathcal{B}_{2}}\Lambda_{1}\tfrac{p}{\rho A}\mathrm{d}z^{1},\qquad u_{\partial,1}^{1}=-y_{c,1}^{\partial,1},\\
u_{c,\partial,1}^{2} & =\int_{\partial\mathcal{B}_{4}}\Lambda_{2}\tfrac{p}{\rho A}\mathrm{d}z^{1},\qquad u_{\partial,1}^{2}=-y_{c,2}^{\partial,1}.
\end{align}
\end{subequations}Hence, if we set the parameters $G_{c,1}^{\partial,1,1}=G_{c,2}^{\partial,1,2}=1$
and $G_{c,2}^{\partial,1,1}=G_{c,1}^{\partial,1,2}=0$, from (\ref{eq:Casimir_cond_BA_diff_op_act_bound_1})
-- (\ref{eq:Casimir_cond_BA_diff_op_unact_bound}) we find
\[
(\Lambda_{1}+\delta_{w}^{\partial,1}\mathcal{C}^{1})|_{\partial\mathcal{B}_{2}}=0,\qquad(\Lambda_{2}+\delta_{w}^{\partial,1}\mathcal{C}^{2})|_{\partial\mathcal{B}_{4}}=0,
\]
whereas $\delta_{w}^{\partial,2}\mathcal{C}^{\lambda}|_{\partial\mathcal{B}_{a}}=0$,
$\delta_{p}^{\partial,2}\mathcal{C}^{\lambda}|_{\partial\mathcal{B}_{a}}=0$
have to be met since $K_{\partial,2}=0$, as well as for the unactuated
boundary
\[
(\dot{w}\delta_{w}^{\partial,1}\mathcal{C}^{\lambda}+\dot{p}\delta_{p}^{\partial,1}\mathcal{C}^{\lambda}+\dot{w}_{[01]}\delta_{w}^{\partial,2}\mathcal{C}^{\lambda}+\dot{p}_{[01]}\delta_{p}^{\partial,2}\mathcal{C}^{\lambda})|_{\partial\mathcal{B}_{u}}=0
\]
for $\lambda=1,2$. Moreover, for the domain condition (\ref{eq:Casimir_cond_BA_diff_op_domain})
we have $\delta_{w}\mathcal{C}^{\lambda}=0$ and $\delta_{p}\mathcal{C}^{\lambda}=0$,
which are trivially satisfied if $\mathcal{C}^{1}$ and $\mathcal{C}^{2}$
stem from total derivatives. In light of this aspects, we find that
a possible choice for Casimir functions is given by $\mathcal{C}^{1}=d_{[01]}(\frac{L_{2}-z^{2}}{L_{2}}\Lambda_{1}w)$
and $\mathcal{C}^{2}=d_{[01]}(\frac{z^{2}}{L_{2}}\Lambda_{2}w)$,
which enable to deduce the relations
\begin{align}
\begin{array}{ccccccc}
x_{c}^{1} & \!= & \!\int_{\partial\mathcal{B}_{2}}\Lambda_{1}w\mathrm{d}z^{1}\!+\!\kappa^{1}, & \: & x_{c}^{2}\! & = & \!\int_{\partial\mathcal{B}_{4}}\Lambda_{2}w\mathrm{d}z^{1}\!+\!\kappa^{2},\end{array}\label{eq:BA_Kirchhoff_plate_controller_state}
\end{align}
with $\kappa^{1}$, $\kappa^{2}$ depending on the initial states
of the plant and the controller. Compared to boundary controllers
for $1$-dimensional systems, like in \cite{Rams2017a} for instance,
this is a major difference, as we have weighted system states integrated
over $\partial\mathcal{B}_{a}$. However, the chosen Casimir functionals
only assign a part of the controller dynamics. In fact, two further
controller states shall be exploited to inject damping, and thus,
the controller is described by
\begin{align*}
J_{c}-R_{c} & =\left[\begin{array}{cccc}
0 & 0 & 0 & 0\\
0 & 0 & 0 & 0\\
0 & 0 & -R_{c}^{33} & J_{c}^{34}-R_{c}^{34}\\
0 & 0 & -J_{c}^{34}-R_{c}^{34} & -R_{c}^{44}
\end{array}\right]\\
G_{c} & =\left[\begin{array}{cc}
1 & 0\\
0 & 1\\
G_{c,1}^{\partial,1,3} & G_{c,2}^{\partial,1,3}\\
G_{c,1}^{\partial,1,4} & G_{c,2}^{\partial,1,4}
\end{array}\right].
\end{align*}
Moreover, to properly shape the closed-loop Hamiltonian $\mathscr{H}_{cl}$,
we set the controller Hamiltonian to
\begin{multline*}
H_{c}=\tfrac{c_{1}}{2}(x_{c}^{1}-x_{c}^{1,d}-\tfrac{u_{s,1}}{c_{1}})^{2}\\
+\tfrac{c_{2}}{2}(x_{c}^{2}-x_{c}^{2,d}-\tfrac{u_{s,2}}{c_{2}})^{2}+\frac{1}{2}M_{c,\mu_{c}\nu_{c}}x_{c}^{\mu_{c}}x_{c}^{\nu_{c}},
\end{multline*}
where $x_{c}^{1,d}$ and $x_{c}^{2,d}$ can be determined by substituting
(\ref{eq:BA_Kirchhoff_Plate_des_equilibrium}) in (\ref{eq:BA_Kirchhoff_plate_controller_state})
and $[M_{c}]$, with $M_{c,\mu_{c}\nu_{c}}\in\mathbb{R}$ for $\mu_{c},\nu_{c}=3,4$,
is a positive definite matrix. Furthermore, we incorporated appropriate
terms yielding constant voltages $u_{s,1}$ and $u_{s,2}$, which
is necessary as (\ref{eq:BA_Kirchhoff_Plate_des_equilibrium}) is
a configuration that requires non-zero power. Finally, we find that
$\mathscr{H}_{cl}$ evolves along closed-loop solutions as $\dot{\mathscr{H}}_{cl}=-x_{c}^{\mu_{c}}M_{c,\mu_{c}\nu_{c}}R_{c}^{\nu_{c}\rho_{c}}M_{c,\rho_{c}\vartheta_{c}}x_{c}^{\vartheta_{c}}\leq0$,
with $\rho_{c},\vartheta_{c}=3,4$.\end{exmp}

\begin{figure}
\centering
\input{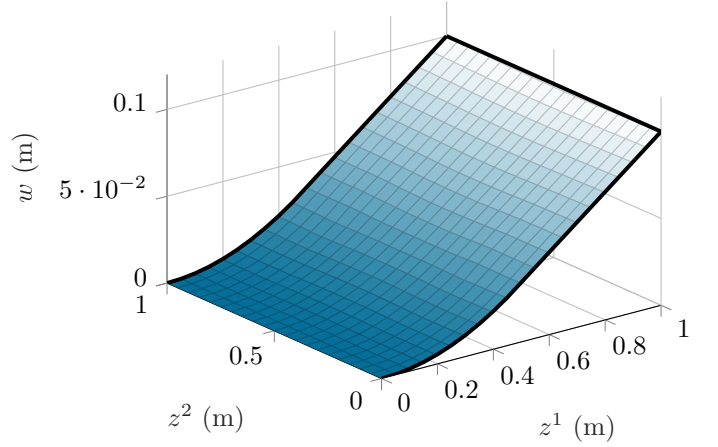}\caption{\label{fig:fin_plate_deflection-1}Final plate deflection $w(z^{1},z^{2},T_{end})$
over $\mathcal{B}$.}
\end{figure}

\begin{figure}
\centering
\input{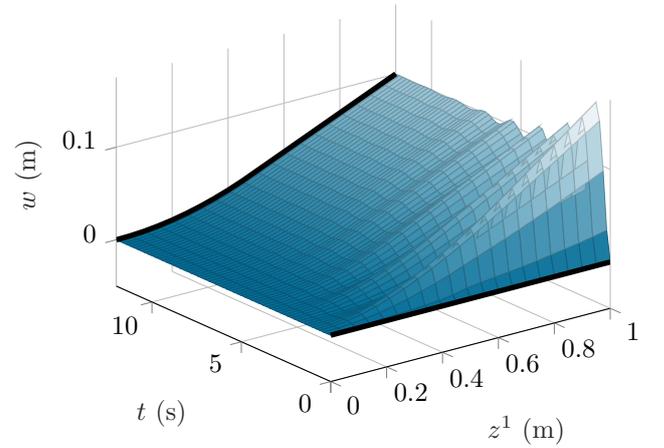}\caption{\label{fig:w_low_bound}Simulation result for the deflection $w$
of the edge $\partial\mathcal{B}_{2}$ over time $t$ and coordinate
$z^{1}$.}
\end{figure}

Let us mention again that detailed stability investigations are not
in the scope of this paper, cf. Rem. \ref{rem:stability}. Thus, we
are content with the achieved findings and with simulation results
presented in Fig. \ref{fig:fin_plate_deflection-1} and \ref{fig:w_low_bound},
which demonstrate the capability of the proposed controller in order
to stabilise the configuration (\ref{eq:BA_Kirchhoff_Plate_des_equilibrium})
with $a=0.1368$ and $b=0.1315$.
 Here, all plate parameters are set to $1$, except Poisson's ratio
that is $\nu=0.2$. Moreover, for the controller parameters we have
chosen $c_{1}=5$, $c_{2}=5$, $J_{c}^{34}=1$, $R_{c}^{34}=1$, $R_{c}^{33}=15$,
$R_{c}^{44}=15$, $G_{c}^{31}=1$, $G_{c}^{42}=1$, $G_{c}^{32}=0$,
$G_{c}^{41}=0$, $M_{c}^{34}=5$, $M_{c}^{33}=25$ $M_{c}^{44}=25$,
$u_{s,1}=-1$ and $u_{s,2}=-1$.

At this point it should be stressed that the controller inputs (\ref{eq:BA_Kirchhoff_Plate_interconnection})
depend on a system state that is distributed over $\partial\mathcal{B}_{2}$
and $\partial\mathcal{B}_{4}$, which cannot be measured. Hence, in
the following the aim is to design an infinite-dimensional observer.

\section{Observer Design}

In this section, we present an energy-based observer-design method,
where the intention is to exploit the pH-system representation and
to introduce error-injection terms based on available measurements
that are located at the boundary. To properly determine these observer-correction
terms, we apply a design approach based on energy balancing -- also
presented in \cite[Sec. 5]{Malzer2018} with regard to the controller
design for spatially $1$-dimensional systems -- for the observer-error
system, which relies on energy shaping and damping injection. For
that purpose, it is assumed to have measurements\footnote{Please note that the two different kinds of measurements denoted by
$\partial,1$ and $\partial,2$ are introduced to take into account
the two different boundary-port categories for systems with $2$nd-order
Hamiltonian.} $\bar{y}_{m}^{\partial,1}$ and $\bar{y}_{m}^{\partial,2}$ (corresponding
to deflections and angular displacements for mechanical systems) that
can be exploited to shape the error Hamiltonian appropriately, and
measurements $y_{m}^{\partial,1}$ and $y_{m}^{\partial,2}$ (velocities
and angular velocities), where the observer-correction terms shall
be introduced such that they are collocated to these measurements
allowing to inject damping into the observer-error system. 

As a first step, we introduce the dynamics of the observer as a copy
of the plant according to
\begin{align}
\dot{\hat{x}}^{\hat{\alpha}} & =(\mathcal{J}^{\hat{\alpha}\hat{\beta}}-\mathcal{R}^{\hat{\alpha}\hat{\beta}})\delta_{\hat{\beta}}\hat{\mathcal{H}},\label{eq:BM_Observer_design_ipH_observer}
\end{align}
with $\hat{\alpha},\hat{\beta}=1,\ldots,n$, where the boundary conditions
need to be determined and $\hat{\mathfrak{H}}=\hat{\mathcal{H}}\Omega$
is the copy of the Hamiltonian density depending on the observer states
$\hat{x}$. Due to the fact that we assume to have measurements available
only at the boundary $\partial\mathcal{B}_{m}=\partial\mathcal{B}_{2}\cup\partial\mathcal{B}_{4}$,
we intend to introduce no observer-correction terms at $\partial\mathcal{B}_{um}=\partial\mathcal{B}_{1}\cup\partial\mathcal{B}_{3}$.
Furthermore, as we have no actuation at this part of the boundary,
no power flow takes place at $\partial\mathcal{B}_{um}$, i.e.
\[
\int_{\partial\mathcal{B}_{um}}\dot{\hat{x}}\rfloor\delta^{\partial,1}\hat{\mathfrak{H}}=0,\qquad\int_{\partial\mathcal{B}_{um}}\dot{\hat{x}}_{[01]}\rfloor\delta^{\partial,2}\hat{\mathfrak{H}}=0.
\]
As a consequence, by exploiting the decomposition Theorem \ref{thm:decomposition_theorem},
the formal change of $\hat{\mathscr{H}}=\int_{\mathcal{B}}\hat{\mathcal{H}}\Omega$
reads
\begin{multline*}
\dot{\hat{\mathscr{H}}}=-\int_{\mathcal{B}}\delta_{\hat{\alpha}}(\hat{\mathcal{H}})\mathcal{R}^{\hat{\alpha}\hat{\beta}}\delta_{\hat{\beta}}(\hat{\mathcal{H}})\Omega\\
+\int_{\partial\mathcal{B}_{m}}(\dot{\hat{x}}^{\hat{\alpha}}\delta_{\hat{\alpha}}^{\partial,1}\hat{\mathcal{H}}+\dot{\hat{x}}_{[01]}^{\hat{\alpha}}\delta_{\hat{\alpha}}^{\partial,2}\hat{\mathcal{H}})\bar{\Omega}_{2},
\end{multline*}
and allows to introduce observer-correction terms as\begin{subequations}\label{eq:BM_Observer_Design_boundary_port_assignment_observer}
\begin{align}
(\delta_{\hat{\alpha}}^{\partial,1}\hat{\mathcal{H}})|_{\partial\mathcal{B}_{m}} & =B_{\hat{\alpha}\varsigma}^{\partial,1}u_{\partial,1}^{\varsigma}-\hat{B}_{\hat{\alpha}\hat{\rho}}^{\partial,1}\hat{k}_{\partial,1}^{\hat{\rho}},\label{eq:BM_Observer_Design_input_assignment_observer1}\\
(\delta_{\hat{\alpha}}^{\partial,2}\hat{\mathcal{H}})|_{\partial\mathcal{B}_{m}} & =B_{\hat{\alpha}\mu}^{\partial,2}u_{\partial,2}^{\mu}-\hat{B}_{\hat{\alpha}\hat{\eta}}^{\partial,2}\hat{k}_{\partial,2}^{\hat{\eta}},\label{eq:BM_Observer_Design_input_assignment_observer2}
\end{align}
with $\hat{\rho}=1,\ldots o_{\partial,1}$, $\hat{\eta}=1,\ldots o_{\partial,2}$
depending on the number of available measurements. The observer inputs
(\ref{eq:BM_Observer_Design_input_assignment_observer1}) and (\ref{eq:BM_Observer_Design_input_assignment_observer2})
comprise the inputs of the plant (\ref{eq:ipH_sys_input_parameterisation})
as well as the observer-correction terms $\hat{k}_{\partial,1}^{\hat{\rho}}$,
$\hat{k}_{\partial,2}^{\hat{\eta}}$, where $\hat{B}_{\hat{\alpha}\hat{\rho}}^{\partial,1}$,
$\hat{B}_{\hat{\alpha}\hat{\eta}}^{\partial,2}$ denote the components
of appropriate mappings $\hat{B}^{\partial,1}$, $\hat{B}^{\partial,2}$,
respectively, which have to be determined depending on the spatial
position of the available measurements. In fact, $\hat{B}_{\hat{\alpha}\hat{\rho}}^{\partial,1}$
and $\hat{B}_{\hat{\alpha}\hat{\eta}}^{\partial,2}$ shall be chosen
such that for the resulting observer-error system the observer-correction
terms $\hat{k}_{\partial,1}^{\hat{\rho}}$, $\hat{k}_{\partial,2}^{\hat{\eta}}$
are collocated to the error terms including the available measurements.
As a consequence, we are able to define
\begin{align}
\begin{array}{ccccccc}
\hat{B}_{\hat{\alpha}\hat{\rho}}^{\partial,1}\dot{\hat{x}}^{\hat{\alpha}}|_{\partial\mathcal{B}_{m}} & = & \hat{y}_{m,\hat{\rho}}^{\partial,1}, & \quad & \hat{B}_{\hat{\alpha}\hat{\eta}}^{\partial,2}\dot{\hat{x}}_{[01]}^{\hat{\alpha}}|_{\partial\mathcal{B}_{m}} & = & \hat{y}_{m,\hat{\eta}}^{\partial,2},\end{array}\label{eq:BM_Observer_Design_output_assignment_observer1}
\end{align}
\end{subequations}which correspond to the observer-equivalent of
the measurements $y_{m,\rho}^{\partial,1}$, $y_{m,\eta}^{\partial,2}$.
Next, we study the observer error $\tilde{x}=x-\hat{x}$, where the
dynamics can be deduced by substituting (\ref{eq:boundary_controlled_ipH_non_diff_op_local})
and (\ref{eq:BM_Observer_design_ipH_observer}) in $\dot{\tilde{x}}=\dot{x}-\dot{\hat{x}}$,
and can be formulated as
\begin{equation}
\dot{\tilde{x}}^{\tilde{\alpha}}=(\mathcal{J}^{\tilde{\alpha}\tilde{\beta}}-\mathcal{R}^{\tilde{\alpha}\tilde{\beta}})\delta_{\tilde{\beta}}\tilde{\mathcal{H}},\label{eq:BM_observer_design_observer_error_System}
\end{equation}
since we confine ourselves to linear systems. To determine the boundary-port
relations of the observer-error system, we study the formal change
of $\tilde{\mathscr{H}}=\int_{\mathcal{B}}\tilde{\mathcal{H}}\Omega$,
where $\tilde{\mathcal{H}}$ exhibits the same form as $\mathcal{H}$
but depends on error coordinates $\tilde{x}$, which follows to
\begin{multline}
\dot{\tilde{\mathscr{H}}}=-\int_{\mathcal{B}}\delta_{\tilde{\alpha}}(\tilde{\mathcal{H}})\mathcal{R}^{\tilde{\alpha}\tilde{\beta}}\delta_{\tilde{\beta}}(\tilde{\mathcal{H}})\Omega\\
+\int_{\partial\mathcal{B}_{m}}(\dot{\tilde{x}}^{\tilde{\alpha}}\delta_{\tilde{\alpha}}^{\partial,1}\tilde{\mathcal{H}}+\dot{\tilde{x}}_{[01]}^{\tilde{\alpha}}\delta_{\tilde{\alpha}}^{\partial,2}\tilde{\mathcal{H}})|_{\partial\mathcal{B}_{m}}\bar{\Omega}_{2}.\label{eq:BM_observer_design_H_punkt_error}
\end{multline}
The restriction to linear systems implies $\delta^{\partial,1}\tilde{\mathfrak{H}}=\delta^{\partial,1}\mathfrak{H}-\delta^{\partial,1}\hat{\mathfrak{H}}$,
$\delta^{\partial,2}\tilde{\mathfrak{H}}=\delta^{\partial,2}\mathfrak{H}-\delta^{\partial,2}\hat{\mathfrak{H}}$,
and thus, as we set the observer inputs to (\ref{eq:BM_Observer_Design_input_assignment_observer1})
and (\ref{eq:BM_Observer_Design_input_assignment_observer2}), regarding
the observer-error system the plant inputs are cancelled. Therefore,
the boundary-inputs read\begin{subequations}\label{eq:BM_Observer_design_inputs_outputs_observer_error}
\begin{equation}
\begin{array}{ccccccc}
(\delta_{\tilde{\alpha}}^{\partial,1}\tilde{\mathcal{H}})|_{\partial\mathcal{B}_{m}} & = & \hat{B}_{\tilde{\alpha}\tilde{\rho}}^{\partial,1}\hat{k}_{\partial,1}^{\tilde{\rho}}, & \; & (\delta_{\tilde{\alpha}}^{\partial,2}\tilde{\mathcal{H}})|_{\partial\mathcal{B}_{m}} & = & \hat{B}_{\tilde{\alpha}\tilde{\eta}}^{\partial,2}\hat{k}_{\partial,2}^{\tilde{\eta}}\end{array},\label{eq:BM_Observer_design_input_assignment_observer_error}
\end{equation}
while the collocated boundary-outputs are given by
\begin{equation}
\begin{array}{ccccccc}
\hat{B}_{\tilde{\alpha}\tilde{\rho}}^{\partial,1}\dot{\tilde{x}}^{\tilde{\alpha}}|_{\partial\mathcal{B}_{m}} & = & \tilde{y}_{\tilde{\rho}}^{\partial,1}, & \quad & \hat{B}_{\tilde{\alpha}\tilde{\eta}}^{\partial,2}\dot{\tilde{x}}_{[01]}^{\tilde{\alpha}}|_{\partial\mathcal{B}_{m}} & = & \tilde{y}_{\tilde{\eta}}^{\partial,2},\end{array}\label{eq:BM_Observer_design_output_assignment_observer_error}
\end{equation}
\end{subequations}where we have $\tilde{y}_{\tilde{\rho}}^{\partial,1}=y_{m,\tilde{\rho}}^{\partial,1}-\hat{y}_{\tilde{\rho}}^{\partial,1}$
and $\tilde{y}_{\tilde{\eta}}^{\partial,2}=y_{m,\tilde{\eta}}^{\partial,2}-\hat{y}_{\tilde{\eta}}^{\partial,2}$
with the measurements $y_{m,\tilde{\rho}}^{\partial,1}$, $y_{m,\tilde{\eta}}^{\partial,2}$
as well as the corresponding observer quantities $\hat{y}_{m,\tilde{\rho}}^{\partial,1}$,
$\hat{y}_{m,\tilde{\eta}}^{\partial,2}$ according to (\ref{eq:BM_Observer_Design_output_assignment_observer1}).
Let us stress again that the coefficients $\hat{B}_{\tilde{\alpha}\tilde{\rho}}^{\partial,1}$
and $\hat{B}_{\tilde{\alpha}\tilde{\eta}}^{\partial,2}$ take into
account the spatial position of the available measurements.

Consequently, the dynamics of the observer error are reformulated
as a pH-system, where the boundary-correction terms $\hat{k}_{\partial,1}$
and $\hat{k}_{\partial,2}$ shall be determined such that the observer-error
system exhibits a desired behaviour. To this end, we apply the energy-balancing
approach presented in \cite[Sec. 5]{Malzer2018}. In particular, the
observer-correction terms $\hat{k}_{\partial,1}$ and $\hat{k}_{\partial,2}$
shall be used to shape the error Hamiltonian $\tilde{\mathscr{H}}$
and to inject damping into the observer-error system.

\begin{thm}[Observer Design]\label{thm:Observer_Design}Consider
the observer- \newline error system (\ref{eq:BM_observer_design_observer_error_System})
with the boundary-inputs and -outputs (\ref{eq:BM_Observer_design_inputs_outputs_observer_error}),
where the observer-correction terms are splitted according to $\hat{k}_{\partial,1}^{\tilde{\rho}}=\hat{\beta}_{\partial,1}^{\tilde{\rho}}+\hat{\gamma}_{\partial,1}^{\tilde{\rho}}$
and $\hat{k}_{\partial,2}^{\tilde{\eta}}=\hat{\beta}_{\partial,2}^{\tilde{\eta}}+\hat{\gamma}_{\partial,2}^{\tilde{\eta}}$.
Thus, if we find a $\tilde{\mathcal{H}}_{a}$ such that the matching
conditions\begin{subequations}\label{eq:BM_Observer_Design_Matching_Conditions}
\begin{align}
(\mathcal{J}^{\tilde{\alpha}\tilde{\beta}}-\mathcal{R}^{\tilde{\alpha}\tilde{\beta}})\delta_{\tilde{\beta}}\tilde{\mathcal{H}}_{a} & =0\label{eq:Matching_Condition_Domain}\\
(\delta_{\tilde{\alpha}}^{\partial,1}\tilde{\mathcal{H}}_{a})|_{\partial\mathcal{B}_{um}} & =0\label{eq:Matching_Condition_boundary_1}\\
(\delta_{\tilde{\alpha}}^{\partial,2}\tilde{\mathcal{H}}_{a})|_{\partial\mathcal{B}_{um}} & =0\label{eq:Matching_Condition_boundary_2}
\end{align}
\end{subequations}are fulfilled, the energy-shaping correction terms\begin{subequations}\label{eq:BM_Observer_design_energy_shaping}
\begin{align}
\int_{\partial\mathcal{B}_{m}}\hat{B}_{\tilde{\alpha}\tilde{\rho}}^{\partial,1}\hat{\beta}_{\partial,1}^{\tilde{\rho}}\bar{\Omega}_{2} & =-\int_{\partial\mathcal{B}_{m}}\delta_{\tilde{\alpha}}^{\partial,1}\tilde{\mathcal{H}}_{a}\bar{\Omega}_{2}\\
\int_{\partial\mathcal{B}_{m}}\hat{B}_{\tilde{\alpha}\tilde{\eta}}^{\partial,2}\hat{\beta}_{\partial,2}^{\tilde{\eta}}\bar{\Omega}_{2} & =-\int_{\partial\mathcal{B}_{m}}\delta_{\tilde{\alpha}}^{\partial,2}\tilde{\mathcal{H}}_{a}\bar{\Omega}_{2}
\end{align}
\end{subequations}map the observer-error system (\ref{eq:BM_observer_design_observer_error_System})
into the target system
\begin{equation}
\dot{\tilde{x}}^{\tilde{\alpha}}=(\mathcal{J}^{\tilde{\alpha}\tilde{\beta}}-\mathcal{R}^{\tilde{\alpha}\tilde{\beta}})\delta_{\tilde{\beta}}\tilde{\mathcal{H}}_{d},\label{eq:BM_Observer_Design_target_system}
\end{equation}
ensuring that $\tilde{\mathscr{H}}_{d}=\int_{\mathcal{B}}\tilde{\mathcal{H}}_{d}\Omega$,
with $\tilde{\mathcal{H}}_{d}=\tilde{\mathcal{H}}+\tilde{\mathcal{H}}_{a}$,
exhibits a certain minimum. Thus, if we set the new inputs of the
target system, which can be parameterised as
\begin{equation}
(\delta_{\tilde{\alpha}}^{\partial,1}\tilde{\mathcal{H}}_{d})|_{\partial\mathcal{B}_{m}}\!=\!\hat{B}_{\tilde{\alpha}\tilde{\rho}}^{\partial,1}\hat{\gamma}_{\partial,1}^{\tilde{\rho}},\;(\delta_{\tilde{\alpha}}^{\partial,2}\tilde{\mathcal{H}}_{d})|_{\partial\mathcal{B}_{m}}\!=\!\hat{B}_{\tilde{\alpha}\tilde{\eta}}^{\partial,2}\hat{\gamma}_{\partial,2}^{\tilde{\eta}},\label{eq:BM_Observer_design_damping_injection_input}
\end{equation}
and are referred to as damping-injection inputs, to
\begin{equation}
\hat{\gamma}_{\partial,1}^{\tilde{\rho}}=-\tilde{K}_{\partial,1}^{\tilde{\rho}\tilde{\kappa}}\tilde{y}_{\tilde{\kappa}}^{\partial,1},\qquad\hat{\gamma}_{\partial,2}^{\tilde{\mu}}=-\tilde{K}_{\partial,2}^{\tilde{\mu}\tilde{\sigma}}\tilde{y}_{\tilde{\sigma}}^{\partial,2},\label{eq:BM_Observer_design_damping_injection_law}
\end{equation}
with appropriate positive definite mappings $\tilde{K}_{\partial,1}$
and $\tilde{K}_{\partial,2}$, the desired error system (\ref{eq:BM_Observer_Design_target_system})
is dissipative.\end{thm}\begin{pf}The matching condition (\ref{eq:Matching_Condition_Domain})
follows immediately by substituting the ansatz $\tilde{\mathcal{H}}_{d}=\tilde{\mathcal{H}}+\tilde{\mathcal{H}}_{a}$
in
\[
(\mathcal{J}^{\tilde{\alpha}\tilde{\beta}}-\mathcal{R}^{\tilde{\alpha}\tilde{\beta}})\delta_{\tilde{\beta}}\tilde{\mathcal{H}}=(\mathcal{J}^{\tilde{\alpha}\tilde{\beta}}-\mathcal{R}^{\tilde{\alpha}\tilde{\beta}})\delta_{\tilde{\beta}}\tilde{\mathcal{H}}_{d},
\]
where (\ref{eq:Matching_Condition_Domain}) is trivially satisfied
if $\tilde{\mathcal{H}}_{a}$ stems from total derivatives. Next,
we compare the power-balance relation of the error-Hamiltonian, which
reads
\begin{multline*}
\dot{\tilde{\mathscr{H}}}=-\int_{\mathcal{B}}\delta_{\tilde{\alpha}}(\tilde{\mathcal{H}})\mathcal{R}^{\tilde{\alpha}\tilde{\beta}}\delta_{\tilde{\beta}}(\tilde{\mathcal{H}})\Omega\\
+\int_{\partial\mathcal{B}_{m}}(\dot{\tilde{x}}^{\tilde{\alpha}}\hat{B}_{\tilde{\alpha}\tilde{\rho}}^{\partial,1}\hat{k}_{\partial,1}^{\tilde{\rho}}+\dot{\tilde{x}}_{[01]}^{\tilde{\alpha}}\hat{B}_{\tilde{\alpha}\tilde{\eta}}^{\partial,2}\hat{k}_{\partial,2}^{\tilde{\eta}})\bar{\Omega}_{2}
\end{multline*}
by means of (\ref{eq:BM_Observer_design_input_assignment_observer_error}),
to the formal change of the desired Hamiltonian $\tilde{\mathscr{H}}_{d}$.
Here, since we cannot shape the Hamiltonian at the boundary $\partial\mathcal{B}_{um}$,
it follows that
\[
\int_{\partial\mathcal{B}_{um}}(\dot{\tilde{x}}^{\tilde{\alpha}}\delta_{\tilde{\alpha}}^{\partial,1}\tilde{\mathcal{H}}_{d}+\dot{\tilde{x}}_{[01]}^{\tilde{\alpha}}\delta_{\tilde{\alpha}}^{\partial,2}\tilde{\mathcal{H}}_{d})\bar{\Omega}_{2}=0,
\]
which allows to find the matching conditions (\ref{eq:Matching_Condition_boundary_1})
and (\ref{eq:Matching_Condition_boundary_2}). Thus, for the formal
change of $\tilde{\mathscr{H}}_{d}$ we have
\begin{multline}
\dot{\tilde{\mathscr{H}}}_{d}=-\int_{\mathcal{B}}\delta_{\tilde{\alpha}}(\tilde{\mathcal{H}}_{d})\mathcal{R}^{\tilde{\alpha}\tilde{\beta}}\delta_{\tilde{\beta}}(\tilde{\mathcal{H}}_{d})\Omega\\
+\int_{\partial\mathcal{B}_{m}}(\dot{\tilde{x}}^{\tilde{\alpha}}\underset{\hat{B}_{\tilde{\alpha}\tilde{\rho}}^{\partial,1}\hat{\gamma}_{\partial,1}^{\tilde{\rho}}}{\underbrace{\delta_{\tilde{\alpha}}^{\partial,1}\tilde{\mathcal{H}}_{d}}}+\dot{\tilde{x}}_{[01]}^{\tilde{\alpha}}\underset{\hat{B}_{\tilde{\alpha}\tilde{\eta}}^{\partial,2}\hat{\gamma}_{\partial,2}^{\tilde{\eta}}}{\underbrace{\delta_{\tilde{\alpha}}^{\partial,2}\tilde{\mathcal{H}}_{d}}})\bar{\Omega}_{2},\label{eq:EBC_H_d_punkt}
\end{multline}
enabling to introduce the inputs according to (\ref{eq:BM_Observer_design_damping_injection_input}).
Therefore, a comparison of the target-system inputs with (\ref{eq:BM_Observer_design_input_assignment_observer_error}),
where we insert $\hat{k}_{\partial,1}^{\tilde{\rho}}=\hat{\beta}_{\partial,1}^{\tilde{\rho}}+\hat{\gamma}_{\partial,1}^{\tilde{\rho}}$
and $\hat{k}_{\partial,2}^{\tilde{\eta}}=\hat{\beta}_{\partial,2}^{\tilde{\eta}}+\hat{\gamma}_{\partial,2}^{\tilde{\eta}}$
allowing to write
\begin{align*}
\hat{B}_{\tilde{\alpha}\tilde{\rho}}^{\partial,1}\hat{\gamma}_{\partial,1}^{\tilde{\rho}} & =(\delta_{\tilde{\alpha}}^{\partial,1}\tilde{\mathcal{H}})|_{\partial\mathcal{B}_{m}}-\hat{B}_{\tilde{\alpha}\tilde{\rho}}^{\partial,1}\hat{\beta}_{\partial,1}^{\tilde{\rho}}\\
\hat{B}_{\tilde{\alpha}\tilde{\eta}}^{\partial,2}\hat{\gamma}_{\partial,2}^{\tilde{\eta}} & =(\delta_{\tilde{\alpha}}^{\partial,2}\tilde{\mathcal{H}})|_{\partial\mathcal{B}_{m}}-\hat{B}_{\tilde{\alpha}\tilde{\eta}}^{\partial,2}\hat{\beta}_{\partial,2}^{\tilde{\eta}},
\end{align*}
yields the energy-shaping correction terms $\hat{B}_{\tilde{\alpha}\tilde{\rho}}^{\partial,1}\hat{\beta}_{\partial,1}^{\tilde{\rho}}=-(\delta_{\tilde{\alpha}}^{\partial,1}\tilde{\mathcal{H}}_{a})|_{\partial\mathcal{B}_{m}}$
and $\hat{B}_{\tilde{\alpha}\tilde{\eta}}^{\partial,2}\hat{\beta}_{\partial,2}^{\tilde{\eta}}=-(\delta_{\tilde{\alpha}}^{\partial,2}\tilde{\mathcal{H}}_{a})|_{\partial\mathcal{B}_{m}}$.
Hence, if we exploit the boundary-input and -output parameterisation
according to (\ref{eq:BM_Observer_design_damping_injection_input})
and (\ref{eq:BM_Observer_design_output_assignment_observer_error}),
respectively, as well as the damping-injection laws (\ref{eq:BM_Observer_design_damping_injection_law}),
relation (\ref{eq:EBC_H_d_punkt}) reads
\begin{multline*}
\dot{\tilde{\mathscr{H}}}_{d}=-\int_{\mathcal{B}}\delta_{\tilde{\alpha}}(\tilde{\mathcal{H}}_{d})\mathcal{R}^{\tilde{\alpha}\tilde{\beta}}\delta_{\tilde{\beta}}(\tilde{\mathcal{H}}_{d})\Omega\\
-\int_{\partial\mathcal{B}_{m}}(y_{\tilde{\varsigma}}^{\partial,1}K_{\partial,1}^{\tilde{\varsigma}\tilde{\kappa}}y_{\tilde{\kappa}}^{\partial,1}+y_{\tilde{\mu}}^{\partial,2}K_{\partial,2}^{\tilde{\mu}\tilde{\sigma}}y_{\tilde{\sigma}}^{\partial,2})\bar{\Omega}_{2}\leq0.
\end{multline*}
\end{pf}

It should be stressed that Theorem \ref{thm:Observer_Design} only
provides a procedure to properly design the observer-correction term,
where the result depends on the choice for $\tilde{\mathcal{H}}_{a}=\sum_{k=1}^{K}\tilde{\mathcal{H}}_{a}^{k}$.
Thus, in the following we demonstrate the observer-design scheme by
means of the Kirchhoff-Love plate of Ex. \ref{ex:Kirchhoff_Love_plate}.

\begin{exmp}[Observer Design for Ex. 2]\label{Ex:Observer_Design_Kirchhoff_Love_plate}Now,
regarding the observer design we assume that the plate deflection
$w(\frac{3L_{1}}{4},0,t)=\bar{y}_{m,1}^{\partial,1}$ together with
the corresponding velocity $\dot{w}(\frac{3L_{1}}{4},0,t)=y_{m,1}^{\partial,1}$
as well as $w(\frac{3L_{1}}{4},L_{2},t)=\bar{y}_{m,2}^{\partial,1}$
and $\dot{w}(\frac{3L_{1}}{4},L_{2},t)=y_{m,2}^{\partial,1}$ are
available as measurement quantities, where the positions are marked
by $\times$ in Fig. \ref{fig:Schematic_Kirchhoff_Love_plate}. Following
the presented approach, we introduce the dynamics of the ipH-observer
as a copy of the plant
\[
\left[\begin{array}{c}
\dot{\hat{w}}\\
\dot{\hat{p}}
\end{array}\right]\!=\!\left[\begin{array}{cc}
0 & 1\\
-1 & 0
\end{array}\right]\left[\begin{array}{c}
\delta_{\hat{w}}\hat{\mathcal{H}}\\
\delta_{\hat{p}}\hat{\mathcal{H}}
\end{array}\right]\!=\!\left[\begin{array}{c}
\tfrac{\hat{p}}{\rho A}\\
-D_{E}(\hat{w}_{[40]}+2\hat{w}_{[22]}+\hat{w}_{[04]})
\end{array}\right],
\]
with the observer density
\begin{multline}
\hat{\mathcal{H}}=\frac{1}{2\rho A}\hat{p}^{2}+\tfrac{1}{2}D_{E}((\hat{w}_{\left[20\right]})^{2}+(\hat{w}_{\left[02\right]})^{2})\\
+\tfrac{1}{2}D_{E}(2\nu\hat{w}_{\left[20\right]}\hat{w}_{\left[02\right]}+2(1-\nu)(\hat{w}_{\left[11\right]})^{2})\label{eq:H_dach_plate}
\end{multline}
depending on the observer states $\hat{w}$ and $\hat{p}$. Note that
in accordance with the plate configuration the observer dynamics are
restricted to the boundary conditions $\hat{w}=0$, $\hat{w}_{[10]}=0$
for $\partial\mathcal{B}_{1}$, $\hat{M}_{1}=0$ for $\partial\mathcal{B}_{2}$,
$\partial\mathcal{B}_{4}$ and $\hat{Q}_{2}=0$, $\hat{M}_{2}=0$
for $\partial\mathcal{B}_{3}$, where the relations $\hat{Q}_{1}=D_{E}(\hat{w}_{[03]}+(2-\nu)\hat{w}_{[21]})$,
$\hat{M}_{1}=-D_{E}(\hat{w}_{[02]}+\nu\hat{w}_{[20]})$, $\hat{Q}_{2}=-D_{E}(\hat{w}_{[30]}+(2-\nu)\hat{w}_{[12]})$
and $\hat{M}_{2}=D_{E}(\hat{w}_{[20]}+\nu\hat{w}_{[02]})$ can be
deduced by applying the boundary operators (\ref{eq:boundary_operators})
to the observer density (\ref{eq:H_dach_plate}). Thus, in light of
(\ref{eq:BM_Observer_Design_boundary_port_assignment_observer}) we
write
\[
\begin{array}{cclcccc}
\hat{Q}_{1}|_{\partial\mathcal{B}_{2}} & = & B_{11}^{\partial,1}u_{\partial,1}^{1}-\hat{B}_{11}^{\partial,1}\hat{k}_{\partial,1}^{1}, & \quad & \hat{B}_{11}^{\partial,1}\dot{\hat{w}}|_{\partial\mathcal{B}_{2}} & = & \hat{y}_{m,1}^{\partial,1}\\
\hat{Q}_{1}|_{\partial\mathcal{B}_{4}} & = & B_{12}^{\partial,1}u_{\partial,1}^{2}-\hat{B}_{12}^{\partial,1}\hat{k}_{\partial,1}^{2}, & \quad & \hat{B}_{12}^{\partial,1}\dot{\hat{w}}|_{\partial\mathcal{B}_{4}} & = & \hat{y}_{m,2}^{\partial,1}
\end{array}
\]
with the components $\hat{B}_{11}^{\partial,1}$, $\hat{B}_{12}^{\partial,1}$
of the mapping $\hat{B}^{\partial,1}$ and the error-injection terms
$\hat{k}_{\partial,1}^{1}$, $\hat{k}_{\partial,1}^{2}$ to be determined.
As we have measurements available at $(z^{1}\!=\!\frac{3L_{1}}{4},z^{2}\!=\!0)$
and $(z^{1}\!=\!\frac{3L_{1}}{4},z^{2}\!=\!L_{2})$, we set $\hat{B}_{11}^{\partial,1}\!=\!\delta(z^{1}\!-\!\frac{3L_{1}}{4})|_{z^{2}=0}$
and $\hat{B}_{12}^{\partial,1}\!=\!\delta(z^{1}\!-\!\frac{3L_{1}}{4})|_{z^{2}=L_{2}}$
with $\delta(\cdot)$ denoting the Dirac delta function. To derive
proper error-injection terms $\hat{k}_{\partial,1}^{1}$, $\hat{k}_{\partial,1}^{2}$,
we study the dynamics of the observer error, which can be introduced
by means of $\tilde{w}=w-\hat{w}$, $\tilde{p}=p-\hat{p}$ together
with the corresponding derivatives and can be written as
\[
\left[\begin{array}{c}
\dot{\tilde{w}}\\
\dot{\tilde{p}}
\end{array}\right]\!=\!\left[\begin{array}{cc}
0 & 1\\
-1 & 0
\end{array}\right]\left[\begin{array}{c}
\delta_{\tilde{w}}\tilde{\mathcal{H}}\\
\delta_{\tilde{p}}\tilde{\mathcal{H}}
\end{array}\right]\!=\!\left[\begin{array}{c}
\tfrac{\tilde{p}}{\rho A}\\
-D_{E}(\tilde{w}_{[40]}+2\tilde{w}_{[22]}+\tilde{w}_{[04]})
\end{array}\right].
\]
Thus, by means of the formal change of the error Hamiltonian $\tilde{\mathscr{H}}$,
which reads $\dot{\tilde{\mathscr{H}}}=\int_{\partial\mathcal{B}_{2}}\tilde{Q}_{1}\dot{\tilde{w}}\mathrm{d}z^{1}+\int_{\partial\mathcal{B}_{4}}\tilde{Q}_{1}\dot{\tilde{w}}\mathrm{d}z^{1}$
with $\tilde{Q}_{1}|_{\partial\mathcal{B}_{2}}=\delta_{w}^{\partial,1}\tilde{\mathcal{H}}|_{\partial\mathcal{B}_{2}}$
and $\tilde{Q}_{1}|_{\partial\mathcal{B}_{4}}=\delta_{w}^{\partial,1}\tilde{\mathcal{H}}|_{\partial\mathcal{B}_{4}}$,
we are able to introduce the boundary-port relations
\[
\begin{array}{ccccccc}
\tilde{Q}_{1}|_{\partial\mathcal{B}_{2}} & = & \hat{B}_{11}^{\partial,1}\hat{k}_{\partial,1}^{1}, & \quad & \hat{B}_{11}^{\partial,1}\dot{\tilde{w}}|_{\partial\mathcal{B}_{2}} & = & \tilde{y}_{1}^{\partial,1}\\
\tilde{Q}_{1}|_{\partial\mathcal{B}_{4}} & = & \hat{B}_{12}^{\partial,1}\hat{k}_{\partial,1}^{2}, & \quad & \hat{B}_{12}^{\partial,1}\dot{\tilde{w}}|_{\partial\mathcal{B}_{4}} & = & \tilde{y}_{2}^{\partial,1}
\end{array}.
\]
In light of the available measurements $w(\frac{3L_{1}}{4},0)$ and
$w(\frac{3L_{1}}{4},L_{2})$, for the energy-balancing scheme we choose\begin{subequations}\label{eq:H_a_observer_error_plate}
\begin{align}
-\tilde{\mathcal{H}}_{a}^{1} & =d_{[01]}(\tfrac{k_{1}(L_{2}-z^{2})}{2L_{2}}\delta(z^{1}-\tfrac{3L_{1}}{4})\tilde{w}^{2}),\\
-\tilde{\mathcal{H}}_{a}^{2} & =d_{[01]}(\tfrac{k_{2}z^{2}}{2L_{2}}\delta(z^{1}-\tfrac{3L_{1}}{4})\tilde{w}^{2})
\end{align}
\end{subequations}where we intentionally write $-\tilde{\mathcal{H}}_{a}$
such that $\int_{\mathcal{B}}\tilde{\mathcal{H}}_{a}^{2}\Omega$ yields
$\int_{\partial\mathcal{B}_{4}}\frac{k_{1}}{2}\delta(z^{1}-\frac{3L_{1}}{4})\tilde{w}^{2}\mathrm{d}z^{1}$
by means of Stoke's Theorem as $\partial_{[01]}\rfloor\Omega=-\mathrm{d}z^{1}$
for instance. Thus, the ansatz (\ref{eq:H_a_observer_error_plate})
fulfils the matching conditions (\ref{eq:BM_Observer_Design_Matching_Conditions})
and yields $\hat{\beta}_{\partial,1}^{1}=-k_{1}\tilde{w}(\frac{3L_{1}}{4},0)$
and $\hat{\beta}_{\partial,1}^{2}=-k_{2}\tilde{w}(\frac{3L_{1}}{4},L_{2})$
by evaluating (\ref{eq:BM_Observer_design_energy_shaping}). Moreover,
if we use the damping-injection laws
\begin{align*}
\hat{\gamma}_{\partial,1}^{1}\! & =\!-\!\tilde{K}_{\partial,1}^{11}\dot{\tilde{w}}(\tfrac{3L_{1}}{4},0)\!=\!-\!\tilde{K}_{\partial,1}^{11}(\dot{w}(\tfrac{3L_{1}}{4},0)\!-\!\dot{\hat{w}}(\tfrac{3L_{1}}{4},0))\\
\hat{\gamma}_{\partial,1}^{2}\! & =\!-\!\tilde{K}_{\partial,1}^{22}\dot{\tilde{w}}(\tfrac{3L_{1}}{4},L_{2})\!=\!-\!\tilde{K}_{\partial,1}^{22}(\dot{w}(\tfrac{3L_{1}}{4},L_{2})\!-\!\dot{\hat{w}}(\tfrac{3L_{1}}{4},L_{2}))
\end{align*}
with $K_{\partial,1}^{11},K_{\partial,1}^{22}>0$, for the target-system
input, we obtain
\[
\dot{\tilde{\mathscr{H}}}_{d}=-K_{\partial,1}^{11}(\dot{\tilde{w}}(\tfrac{3L_{1}}{4},0))^{2}-K_{\partial,1}^{22}(\dot{\tilde{w}}(\tfrac{3L_{1}}{4},L_{2}))^{2}\leq0,
\]
highlighting that we are able to inject damping into the observer-error
system.

Of course, regarding the observer design a rigorous proof of stability
would be desirable in order to ensure the convergence of the observer.
However, similar to the controller design, cf. Rem. \ref{rem:stability},
here we are content with a non-increasing error-Hamiltonian ensured
by $\dot{\tilde{\mathscr{H}}}_{d}\leq0$ and simulation results. To
this end, the infinite-dimensional observer, where the dynamics
\begin{align*}
\dot{\hat{w}} & =\tfrac{\hat{p}}{\rho A}\\
\dot{\hat{p}} & =-D_{E}(\hat{w}_{[40]}+2\hat{w}_{[22]}+\hat{w}_{[04]})
\end{align*}
are subjected to the boundary conditions
\[
\begin{array}{rcllcc}
\left.\begin{array}{rcc}
\hat{w} & = & 0\\
\hat{w}_{[10]} & = & 0
\end{array}\right\}  & \!\text{for}\! & \partial\mathcal{B}_{1}, & \left.\begin{array}{rcl}
\hat{Q}_{1} & = & \Lambda_{1}u_{in}^{1}-\hat{B}_{11}^{\partial,1}\hat{k}_{\partial,1}^{1}\\
\hat{M}_{1} & = & 0
\end{array}\right\}  & \!\text{for}\! & \partial\mathcal{B}_{2}\\
\left.\begin{array}{rcc}
\hat{Q}_{2} & = & 0\\
\hat{M}_{2} & = & 0
\end{array}\right\}  & \!\text{for}\! & \partial\mathcal{B}_{3}, & \left.\begin{array}{rcl}
\hat{Q}_{1} & = & \Lambda_{2}u_{in}^{2}-\hat{B}_{12}^{\partial,1}\hat{k}_{\partial,1}^{2}\\
\hat{M}_{1} & = & 0
\end{array}\right\}  & \!\text{for}\! & \partial\mathcal{B}_{4},
\end{array}
\]
with the observer-correction terms
\begin{align*}
\hat{k}_{\partial,1}^{1} & \!=\!-\!k_{1}(\bar{y}_{m,1}^{\partial,1}\!-\!\hat{w}(\tfrac{3L_{1}}{4},0))\!-\!\tilde{K}_{\partial,1}^{11}(y_{m,1}^{\partial,1}\!-\!\tfrac{\hat{p}}{\rho A}(\tfrac{3L_{1}}{4},0))\\
\hat{k}_{\partial,1}^{2} & \!=\!-\!k_{2}(\bar{y}_{m,2}^{\partial,1}\!-\!\hat{w}(\tfrac{3L_{1}}{4},L_{2}))\!-\!\tilde{K}_{\partial,1}^{22}(y_{m,2}^{\partial,1}\!-\!\tfrac{\hat{p}}{\rho A}(\tfrac{3L_{1}}{4},L_{2}))
\end{align*}
injected at the boundaries $\partial\mathcal{B}_{2}$ and $\partial\mathcal{B}_{4}$,
respectively, is implemented by means of the finite difference-coefficient
method and initialised with $\hat{w}(z^{1},z^{2},0)=-d(z^{1})^{2}$,
$\hat{p}(z^{1},z^{2},0)=0$ with $d=0.05$. Note that the application
of an observer basically destroys the relation $\dot{\mathscr{C}}^{\lambda}=0$,
i.e. a fast convergence of the observer state is of particular importance,
and therefore, we set $k_{1}=k_{2}=2000$, $\tilde{K}_{\partial,1}^{11}=\tilde{K}_{\partial,1}^{22}=2000$.
The simulation results presented in Fig. \ref{fig:comparison_observer},
where a comparison measurement $w(L_{1},\tfrac{L_{2}}{2})$ -- marked
by $\circ$ in Fig. \ref{fig:Schematic_Kirchhoff_Love_plate} --
and the corresponding observer quantity $\hat{w}(L_{1},\tfrac{L_{2}}{2})$
are depicted, demonstrate the applicability of the combination of
the proposed observer and controller in order to stabilise the configuration
(\ref{eq:BA_Kirchhoff_Plate_des_equilibrium}) depicted in Fig. \ref{fig:fin_plate_deflection-1}.
As the red, dotdashed line in Fig. \ref{fig:comparison_observer}
represents the simulation result $w(L_{1},\tfrac{L_{2}}{2})$ without
using an observer, one find that the combination of observer and controller,
where the inputs of the Casimir-based controller read
\[
\begin{array}{ccccccc}
u_{c,\partial,1}^{1} & = & \int_{\partial\mathcal{B}_{2}}\Lambda_{1}\tfrac{\hat{p}}{\rho A}\mathrm{d}z^{1}, & \quad & u_{c,\partial,1}^{2} & = & \int_{\partial\mathcal{B}_{4}}\Lambda_{2}\tfrac{\hat{p}}{\rho A}\mathrm{d}z^{1}\end{array},
\]
 achieves a similar performance.\end{exmp}
\begin{figure}

\centering
\input{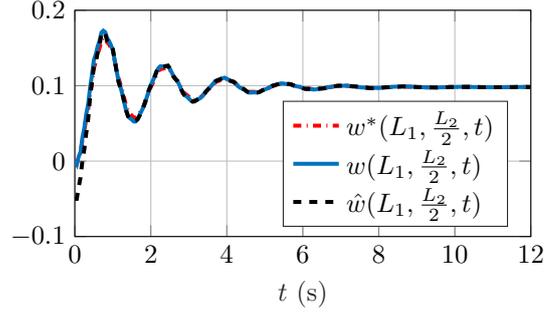}\caption{\label{fig:comparison_observer}Comparison of the measurement $w(L_{1},\frac{L_{2}}{2},t)$
and the observer state $\hat{w}(L_{1},\frac{L_{2}}{2},t)$.}
\end{figure}



\bibliography{my_bib}             

\end{document}